\documentclass[final,1p,times]{elsarticle}

\usepackage[colorlinks=true]{hyperref}
\usepackage{amsfonts}
\usepackage[all]{xy}
\usepackage[mathscr]{eucal}
\usepackage{extarrows}
\usepackage{amsthm,amsmath,amssymb,mathrsfs,amsfonts,graphicx,graphics,latexsym,exscale,cmmib57,dsfont,bbm,amscd}
\usepackage{color,soul}
\newtheorem{theorem}{Theorem}[section]
\newtheorem{proposition}[theorem]{Proposition}%
\newtheorem{remark}[theorem]{Remark}%
\newtheorem{lemma}[theorem]{Lemma}
\newtheorem{definition}[theorem]{Definition}%
\journal{xxx}

\begin{document}
\begin{sloppypar}

\begin{frontmatter}
\title{Variational principle for neutralized packing pressure on subsets}

\author{Zubiao Xiao\corref{cor1}}
\ead{xzb2020@fzu.edu.cn}
\cortext[cor1]{Corresponding author}
\address{School of Mathematics and Statistics, Fuzhou University, Fuzhou 350116, People's Republic of China}
\author{Hongwei Jia}
\ead{jiahongwei2878@163.com}
\address{School of Mathematics and Statistics, Fuzhou University, Fuzhou 350116, People's Republic of China}

\begin{abstract}
In this paper, we introduce the notions of neutralized packing pressures and neutralized measure-theoretic pressures on subsets for a finitely generated  free semigroup action. Let $X$ be a compact metric space and $\mathcal{G}$ be a finite family of continuous self-maps on $X$. We consider the semigroup $G$ generated by $\mathcal{G}$ on $X$. We show that the variational principle between the neutralized  packing pressures $P^{P}_{\mathcal{G}}(Z,f)$ and the neutralized  measure--theoretic upper pressures $\overline{P}_{\mu,{\mathcal{G}} }(Z,f)$ for a given continuous function $f$ and a compact subset $Z \subset X$:

$$P^{P}_{\mathcal{G}}(Z,f)=\lim_{\varepsilon \to 0}\sup \{ \overline{P}_{\mu,\mathcal{G} }(Z,f,\varepsilon):\mu \in M(X), \ \mu(Z)=1 \}.$$
\end{abstract}

\end{frontmatter}

\section{Introduction}\label{sec1}

Topological entropy was first introduced for measuring the complexity of a dynamical system by Adler et al. \cite{Ad} in 1965.  Later, Bowen \cite{Bo} gave equivalent definitions for this entropy when $X$ is a metrizable space and then introduced a definition of topological entropy on subsets. Inspired by Hausdorff dimension \cite{Bowen}, the topological entropy is known as \textit{Bowen topological entropy}. In 1984, Pesin and Pitskel \cite{Pe2} conducted a further research on this kind of entropy. In general, it is difficult to study the properties of dynamical systems for arbitrary group actions, so one needs to find some special group that meet certain conditions, e.g., amenable group, sofic group, etc. The notion of topological entropy for free semigroup actions was introduced by Bi\'{s} \cite{Bi} and Bufetov \cite{Bu}. Later, Wang and Ma \cite{Wang} gave  dimensional characterizations for entropy and generalized Bi\'{s}  topological entropy  to a semigroup action on a metric space. Lin et al. \cite{Lin} introduced a kind of measure-theoretic entropy for the action of a free semigroup and obtained a partial variational principle.

Topological pressure, as a natural extension of topological entropy, makes substantial contributions to dynamical systems. Ruelle \cite{Ru} first introduced the concept of topological pressure of additive potentials for expansive dynamical systems. Walters \cite{Wa} then extended this concept to a compact space with continuous transformation. Furthermore,  Ma and Wu in \cite{Ma1} introduced the topological pressure for the action of a semigroup using the theory of Carath\'{e}odory--Pesin structure. Carvalho et al. \cite{Car} proposed a definition of topological pressure for finitely generated free semigroup actions with a random walk and also provided an associated variational principle. For the progress in topological pressure, the reader can refer to \cite{Cui,Xiao1,Zh2}.
\textit{Packing topological entropy}, introduced by Feng and Huang \cite{Fe}, serves as a dynamical analogy of the packing dimension and a counterpart of Bowen topological entropy. In \cite{Fe}, they also obtained variational principle for packing entropy. Since then, Feng and Huang's variational principles have been extended to different systems and topological pressures. The reader can refer to \cite{Dou,Huang,Ta,WangC,Xiao,Yin,Zheng1,Zh1} for more details. 

Recently, Ben Ovadia and Rodriguez--Hertz \cite{Be} defined the neutralized Bowen open ball for an autonomous dynamical system $(X,f)$ on a compact metric space $X$ with metric $d$, that is to say
$$B_{n}(x,e^{-n\varepsilon})=\{y\in X :d(f^{j}(x),f^{j}(y)<\varepsilon,0\leq j<n)\},$$
where $\varepsilon>0$, $n \in \mathbb{N}$ and $x\in X$. Ben and Rodriguez--Hertz introduced the concept of neutralized Bowen open balls $B_{n}(x,e^{-n\varepsilon})$  as an alternative to the standard Bowen open balls $B_{n}(x,\varepsilon)$. Using these, they defined the neutralized Brin--Katok local entropy, which neutralizes subexponential effects, to estimate the asymptotic measure of sets with a distinctive geometric shape. Furthermore, neutralized Bowen open balls are more useful to describe  the neighborhood with a local linearization of the dynamics. For the relevant knowledge on neutralized entropy and neutralized pressure, readers can refer to \cite{Na,Ya,Sar,Qu}.

In \cite{Zh1}, Zhong and Chen showed that the Bowen topological pressure can be determined by the measure-theoretic pressure. Then Dou et al. \cite{Dou} investigated a variational principle of the packing entropy for the action of a amenable group. Inspired by the above results, we wonder whether the variational principle of the neutralized packing pressure still holds for the action of a free semigroup.

We arrange the rest of this paper as follows. In Section \ref{sec2}, we give the definitions of neutralized packing pressure and neutralized measure--theoretic upper pressure for the action of a free semigroup and study some properties. In Section \ref{sec3},  we define the neutralized Katok’s packing pressure, and give a variational principle between the neutralized packing pressure and the neutralized upper measure--theoretic pressure for the action of a free semigroup.

\section{Neutralized packing pressure for the action of a free semigroup} \label{sec2}

Let $X$ be a compact metric space with metric $d$ and $\mathcal{G}=\left\{f_{1},f_{2},\dots,f_{k} \right\}$ be a finite family of continuous self--maps on $X$. We consider the free semigroup $G$  with the finite generator set $\mathcal{G}$ on $X$, that is, $G=\bigcup_{n \in \mathbb{N}} G_{n}$, where 
$$G_{n}=\left\{g_{1} \circ \cdots \circ g_{m}: 1 \le m < n \text{ and } g_{1},\dots,g_{m} \in \mathcal{G} \right\} \cup \left\{\text{id}_{X}\right\}.$$
We can check that $\left| G_{n} \right|= \sum_{i=0}^{n-1}k^{i}$ for each $n \ge 1$.

For any $n\in \mathbb{N}$, we define
$$d_{n}(x,y)=\max_{g \in G_{n}} d(gx,gy),\ \forall x,y \in X.$$
It is easy to see that $d_{F}$ is a metric on $X$.
For $\varepsilon>0$ and $x \in X$, we denote
$$B_{n}(x,\varepsilon )=\left\{y \in X: d_{n}(x,y) <\varepsilon \right\}$$
and
$$\overline{B}_{n}(x,\varepsilon )=\left\{y \in X: d_{n}(x,y) \leq \varepsilon \right\},$$
which are respectively the \textit{open and closed Bowen balls} with the center $x$ and the radius $\varepsilon$.

Define neutralized Bowen balls as follow:
$$\begin{aligned}
   B_{n}(x,e^{-n\varepsilon} )&=\left\{y \in X: d_{n}(x,y) <e^{-n\varepsilon} \right\},\\
   \overline{B}_{n}(x,e^{-n\varepsilon} )&=\left\{y \in X: d_{n}(x,y) \leq e^{-n\varepsilon} \right\}.
\end{aligned}$$

Let $f \in C(X,\mathbb{R})$, where $C(X,\mathbb{R})$ denote the set of all continuous functions of $X$. For $x \in X$, we write
$$\begin{aligned}
   f_{n}(x)&=\sum_{g \in G_{n}}f(g(x)),\\
   f_{n}(x,\varepsilon )&=\sup_{y \in B_{n}(x,e^{-n\varepsilon} )}f_{n}(y),\\
   \overline{f}_{n}(x,\varepsilon )&=\sup_{y \in \overline{B}_{n}(x,e^{-n\varepsilon} )}f_{n}(y).
\end{aligned}$$
For any $Z\subset X$, $n \in \mathbb{N}$, $\alpha \in \mathbb{R}$, $\varepsilon > 0$, and $f \in C(X,\mathbb{R})$, we define
\begin{equation}
M^{P}_{\mathcal{G}}(n,\alpha ,\varepsilon ,Z,f)= \sup\left \{ \sum_{i} e^{-\alpha |G_{n_{i}}|+f_{n_{i}} (x_{i}) }  \right \}, \label{eq:1.1}
\end{equation}
where the supremum is taken over all finite or countable pairwise disjoint collections of $\left\{\overline{B}_{n_{i}} (x_{i},e^{-n\varepsilon})\right\}_{i}$, such that $x_{i} \in Z$, $n_{i}\geq n$ for all $i$. The quantity $M^{P}_{\mathcal{G}}(n,\alpha ,\varepsilon ,Z,f)$ does not increase as $n$ increases, therefore the following limit exists:
$$M^{P}_{\mathcal{G}}(\alpha ,\varepsilon ,Z,f)=\lim_{n\to +\infty}M^{P}_{\mathcal{G}}(n,\alpha ,\varepsilon ,Z,f).$$

Define
$$M^{\mathcal{P} }_{\mathcal{G}}(\alpha ,\varepsilon ,Z,f)=\inf\left \{ \sum_{i=1}^{\infty } M^{P}_{\mathcal{G}}(\alpha ,\varepsilon ,Z_{i},f) :Z\subseteq \bigcup_{i=1}^{\infty} Z_{i}\right \}.$$
It is easy to check that there exists a critical value of the parameter $\alpha$ , when $\alpha$ goes from $-\infty$ to $+\infty$, such that the quantity $M^{P}_{\mathcal{G}}(\alpha ,\varepsilon ,Z,f)$ jump from $+\infty$ to $0$. Hence, we can define the value

$$\begin{aligned}
P^{P}_{\mathcal{G}}(\varepsilon ,Z,f)
&=\sup\left\{\alpha:M^{\mathcal{P} }_{\mathcal{G}}(\alpha ,\varepsilon ,Z,f)=+\infty\right\}\\
&=\inf\left\{\alpha:M^{\mathcal{P} }_{\mathcal{G}}(\alpha ,\varepsilon ,Z,f)=0\right\}.
\end{aligned}$$
It is not difficult to see that $P^{P}_{\mathcal{G}}(\varepsilon ,Z,f)$ not increases when $\varepsilon$ decreases.
\begin{definition}
   We call the following quantity
   $$
     P^{P}_{\mathcal{G}}(Z,f)=\lim_{\varepsilon \to 0} P^{P}_{\mathcal{G}}(\varepsilon ,Z,f)$$
     a \textit{neutralized packing topological pressure} of $\mathcal{G}$ on the subset $Z$ with respect to $f$.
\end{definition}

Replacing $f_{{n_{i} } } (x_{i})$ in equations \eqref{eq:1.1} by $\overline{f}_{F_{n_{i} } } (x_{i},\varepsilon )$, we can define new functions $\mathcal{M}_{\mathcal{G}}^{P}$. Then we denote the critical value by $P^{P'}_{\mathcal{G}}(\varepsilon ,Z,f).$

we will investigate some basic properties of neutralized packing pressure.

\begin{proposition}\label{p1}
Let $G$ be a free semigroup acting on a compact metric space $(X,d)$ with a finite generator set $\mathcal{G}=\left\{f_{1},f_{2},\dots,f_{k} \right\}$. If $f \in C(X,\mathbb{R})$ and $Z \subseteq X$, then we have
$$P^{P}_{\mathcal{G}}(Z,f)=\lim_{\varepsilon \to 0} P^{P'}_{\mathcal{G}}(\varepsilon ,Z,f).$$
\end{proposition}
\begin{proof}
   It is clear that
   $P^{P}_{\mathcal{G}}(Z,f)\leq \lim_{\varepsilon \to 0} P^{P'}_{\mathcal{G}}(\varepsilon ,Z,f).$
   We shall prove the converse inequality. Let
   $$\gamma (\varepsilon)=\lim_{n \to +\infty}\gamma_{n} (\varepsilon),$$
   where $\gamma_{n} (\varepsilon)=\sup\left\{\left|f(x)-f(y)\right|:x,y \in X,d(x,y)\leq 2e^{-n\varepsilon} \right\}.$
   
   For any $n \in \mathbb{N}$, $u,v \in \overline{B}_{{n_{i}}}(x_{i},e^{-n\varepsilon} )$, and $g \in G_{n_{i}}$, we have $d(g(u),g(v)) \leq 2e^{-n\varepsilon}$, and then
   $$\left|f(g(u))-f(g(v))\right| \leq \gamma_{n}(\varepsilon).$$
   It follows that
   $$|G_{n_{i}}| \gamma_{n}(\varepsilon)\geq \overline{f}_{{n_{i} } } (x_{i},\varepsilon )- f_{{n_{i} } } (x_{i}), \  \forall n_{i}\geq n.$$
   Then we have
   $$\begin{aligned}
   M^{P}_{\mathcal{G}}(n,\alpha ,\varepsilon ,Z,f)&= \sup\left \{ \sum_{i} e^{-\alpha |G_{n_{i}}|+f_{n_{i}} (x_{i}) }  \right \}\\
   &\ge \sup\ \left \{ \sum_{i} e^{-(\alpha+\gamma_{n}(\varepsilon)) |G_{n_{i}}|+\overline{f}_{{n_{i} } } (x_{i},\varepsilon) }\right \}\\
   &=\mathcal{M}^{P}_{\mathcal{G}}(n,\alpha+\gamma_{n}(\varepsilon) ,\varepsilon ,Z,f).
\end{aligned}$$
It follows that $P^{P}_{\mathcal{G}}(\varepsilon,Z,f)\ge P^{P'}_{\mathcal{G}}(\varepsilon ,Z,f)-\gamma(\varepsilon)$. Since $f$ is uniformly continuous on $X$, we know that $\gamma(\varepsilon)=0$.
Then we obtain the desired inequality by letting $\varepsilon \to 0$.
\end{proof}

\begin{proposition}\label{p2}
Let $G$ be a free semigroup acting on a compact metric space $(X,d)$ with a finite generator set $\mathcal{G}=\left\{f_{1},f_{2},\dots,f_{k} \right\}$. If $f \in C(X,\mathbb{R})$ and $Z \subseteq X$, then we have
\begin{enumerate}[(1)]
\item  The value of $P^{P}_{\mathcal{G}}(Z,f)$ is independent of the choice of metrics on $X$. \label{p2.21}
\item if $Z_{1} \subseteq Z_{2}$, then $P^{P}_{\mathcal{G}}(Z_{1},f) \leq P^{P}_{\mathcal{G}}(Z_{2},f)$.\label{p2.2}
\item if $Z=\bigcup_{i \in I}Z_{i}$ with I at most countable, then
\begin{enumerate}[(3-a)]
\item $M^{\mathcal{P}}_{\mathcal{G}}(\alpha ,\varepsilon ,Z,f) \leq \sum_{i \in I}M^{\mathcal{P}}_{\mathcal{G}}(\alpha ,\varepsilon ,Z_{i},f,);$ \label{p2.3a}
\item $P^{P}_{\mathcal{G}}(Z,f)=\sup_{i \in I}P^{P}_{\mathcal{G}}(Z_{i},f).$ \label{p2.3b}
\end{enumerate}
\end{enumerate}
\end{proposition}
\begin{proof}
   (\ref{p2.2}) can be obtained directly by the definitions of the neutralized packing topological pressures. We now show (\ref{p2.21}). 

   Let $d_{1}$ and $d_{2}$ be two compatible metrics on $X$. Then for every $\varepsilon'>0$, there is $\delta'>0$, such that $d_{1}(x,y)<\delta'$ implies $d_{2}(x,y)<\varepsilon'$ for all $x,y \in X$. Now fix $\varepsilon >0$ and let $0< \delta< \varepsilon$. For every $n \in \mathbb{N}$, there is $n_{0} >n$ so that $d_{1}(x,y)< e^{-n_{0}\delta}$ implies $d_{2}(x,y)< e^{-n\varepsilon}$ for all $x,y \in X$. Hence we have $M^{P,{d_{2}}}_{\mathcal{G}}(n,\alpha ,\varepsilon ,Z,f) \leq M^{P,d_{1}}_{\mathcal{G}}(n_{0},\alpha ,\delta,Z,f)$ and so $M^{P,{d_{2}}}_{\mathcal{G}}(\alpha ,\varepsilon ,Z,f) \leq M^{P,d_{1}}_{\mathcal{G}}(\alpha ,\delta,Z,f)$.
   Then we have $P^{P,{d_{2}}}_{\mathcal{G}}(\varepsilon ,Z,f) \leq P^{P,d_{1}}_{\mathcal{G}}(\delta ,Z,f)$. As $\varepsilon \to 0$, one has $P^{P,d_{2}}_{\mathcal{G}}(Z,f) \leq P^{P,d_{1}}_{\mathcal{G}}(Z,f)$. we get the converse inequality by exchanging the role of $d_{1}$ and $d_{2}$.

   We shall prove (3-a). Given $\gamma >0$ and $i \in I $, we can find $\left\{Z_{i,j}\right\}_{j \geq 0}$ such that $Z_{i} \subseteq  \bigcup_{j \geq 0}Z_{i,j}$ and
   $$\sum_{j \geq 0}M^{{P}}_{\mathcal{G}}(\alpha ,\varepsilon ,Z_{i,j},f) \leq M^{\mathcal{P}}_{\mathcal{G}}(\alpha ,\varepsilon ,Z_{i},f)+\frac{\gamma}{2^{i}}.$$ Thus
   $$\sum_{i \in I}\sum_{j \geq 0}M^{{P}}_{\mathcal{G}}(\alpha ,\varepsilon ,Z_{i,j},f) \leq\sum_{i \in I} M^{\mathcal{P}}_{\mathcal{G}}(\alpha ,\varepsilon ,Z_{i},f)+2\gamma.$$
   Then it is obvious that
   $$M^{\mathcal{P}}_{\mathcal{G}}(\alpha ,\varepsilon ,Z,f) \leq \sum_{i \in I}\sum_{j \geq 0}M^{{P}}_{\mathcal{G}}(\alpha ,\varepsilon ,Z_{i,j},f).$$
   Letting $\gamma \to 0$, the desired inequality follows.

   We now show that (3-b) holds. If $\sup_{i \in I}P^{P}_{\mathcal{G}}(Z_{i},f)<s$, then there exists $\varepsilon_{0}>0$, 
   for any $0<\varepsilon<\varepsilon_{0}$ and $i \in I$, $P^{P}_{\mathcal{G}}(\varepsilon,Z_{i},f)<s$, and thus $M^{\mathcal{P}}_{\mathcal{G}}(s,\varepsilon ,Z_{i},f)=0$, which implies $M^{\mathcal{P}}_{\mathcal{G}}(s,\varepsilon ,Z,f)=0$(utilizing (3-a)). Hence $P^{P}_{\mathcal{G}}(\varepsilon,Z,f) \leq s$. It follows that
   $$P^{P}_{\mathcal{G}}(Z,f) \leq \sup_{i \in I}P^{P}_{\mathcal{G}}(Z_{i},f).$$
   The opposite inequality follows from (\ref{p2.2}).
\end{proof}

Next, we will introduce the definition of neutralized measure--theoretic upper pressure for the action of a free semigroup.

Let $M(X)$ denote the set of all Borel probability measures on $X$. For $f \in C(X,\mathbb{R})$ and $\mu \in M(X)$, we define
$$\overline{P}_{\mu ,{\mathcal{G}}}(x,f,\varepsilon):= \limsup_{n \to +\infty} \frac{-\log{\mu}(B_{n}(x,e^{-n\varepsilon} )) +f_{n}(x)}{| G_{n} | }$$
and
$$\overline{P}_{\mu ,{\mathcal{G}}}(Z,f,\varepsilon):=\int_{Z} \overline{P}_{\mu,{\mathcal{G}} }(x,f,\varepsilon) \ \text{d}\mu (x).$$
\begin{definition}\label{d1}
The \textit{neutralized measure--theoretic upper pressures} on the subset $Z$ of $\mu$ for ${\mathcal{G}}$ with respect to $f$ is defined as
$$\overline{P}_{\mu,{\mathcal{G}} }(Z,f)= \lim_{\varepsilon  \to 0}\overline{P}_{\mu,{\mathcal{G}} }(Z,f,\varepsilon).$$
\end{definition}

\begin{remark}
In fact, for every $\alpha>0$, the set $\{x\in X : \mu (B_{n}(x,e^{-n\varepsilon}))> \alpha \}$ is open. This implies that the mappings $x \mapsto{ \mu (B_{n}(x,e^{n\varepsilon}))}$ and $x \mapsto \overline{P}_{\mu ,{\mathcal{G}}}(x,f,\varepsilon)$ are measurable, therefore $\overline{P}_{\mu ,{\mathcal{G}}}(x,f,\varepsilon)$  is integrable and $\overline{P}_{\mu ,{\mathcal{G}}}(Z,f,\varepsilon)$ is well-defined.
\end{remark}

\section{Variational Principle}\label{sec3}

In this section,we will prove the variational principle among the neutralized packing topological pressure and the neutralized measure-theoretic upper pressures.
\begin{theorem}\label{T3}
Let $G$ be a free semigroup acting on a compact metric space $(X,d)$ with a finite generator set $\mathcal{G}=\left\{f_{1},f_{2},\dots,f_{k} \right\}$. If $f \in C(X,\mathbb{R})$, $Z \subset X$ is a nonempty compact set and $ P^{P}_{\mathcal{G}}(Z,f) > \|f\|$, where $\|f\|:=\sup_{x \in X}|f(x)|$, then
$$P^{P}_{\mathcal{G}}(Z,f)=\lim_{\varepsilon \to 0}\sup \{ \overline{P}_{\mu,\mathcal{G} }(Z,f,\varepsilon):\mu \in M(X), \ \mu(Z)=1 \}.$$
\end{theorem}

\subsection{Lower bound}\label{subsec2}

The following theorem is the lower bound part of the variational principle.
\begin{theorem}\label{T4}
Let $G$ be a free semigroup acting on a compact metric space $(X,d)$ with a finite generator set $\mathcal{G}=\left\{f_{1},f_{2},\dots,f_{k} \right\}$. If $Z\subseteq X$ is a Borel set and $f \in C(X,\mathbb{R})$. Then we have
$$P^{P}_{\mathcal{G}}(Z,f)\geq \lim_{\varepsilon \to 0}\sup \{ \overline{P}_{\mu,\mathcal{G} }(Z,f,\varepsilon):\mu \in M(X), \ \mu(Z)=1 \}.$$
\end{theorem}

To obtain a variational principle, we need to define the neutralized Katok’s packing pressure  of Borel probability measure. A parallel methodology is applicable to the classical Katok's entropy which is defined via spanning sets \cite{Katok}.

Given $\mu \in M(X)$,  $\alpha \in \mathbb{R}$, $\varepsilon>0$, $0<\delta<1$ and $f \in C(X,\mathbb{R})$, we define
$$M_{\mu,\mathcal{G}}^{\mathcal{P}}(\alpha ,\varepsilon ,\delta,f)=\inf\ \left \{\sum_{i=1}^{\infty } M^{P}_{\mathcal{G}}(\alpha ,\varepsilon ,Z_{i},f) :\mu(\bigcup_{i=1}^{\infty}Z_{i}) \geq 1-\delta \right \},$$
when $\alpha$ goes from $-\infty$ to $+\infty$,
$M_{\mathcal{G},\mu}^{\mathcal{P}}(\alpha ,\varepsilon ,\delta,f)$
jump from $+\infty$ to $0$ at critical values.
Hence we can define the numbers
$$\begin{aligned} P^{KP}_{\mu,\mathcal{G}}(\varepsilon ,\delta,f)
   &=\sup\left\{\alpha:M_{\mu,\mathcal{G}}^{\mathcal{P}}(\alpha ,\varepsilon ,\delta,f)=+\infty\right\}\\
   &=\inf\left\{\alpha:M_{\mu,\mathcal{G}}^{\mathcal{P}}(\alpha ,\varepsilon ,\delta,f)=0\right\}.
\end{aligned}$$
Let $P^{KP}_{\mu,\mathcal{G}}(\varepsilon,f)=\lim_{\delta \to 0}P^{KP}_{\mu,\mathcal{G}}(\varepsilon ,\delta,f).$
\begin{definition}
We call the following quantity
$$P^{KP}_{\mu,\mathcal{G}}(f)=\lim_{\varepsilon \to 0}P^{KP}_{\mu,\mathcal{G}}(\varepsilon,f)$$
a \textit{neutralized Katok’s packing pressure} of $\mu$ for $\mathcal{G}$ with respect to $f$.
\end{definition}
We shall show that the neutralized measure--theoretic upper pressures is a lower bound of the neutralized Katok’s packing pressure.

\begin{proposition}\label{p3}
Let $G$ be a free semigroup acting on a compact metric space $(X,d)$ with a finite generator set $\mathcal{G}=\left\{f_{1},f_{2},\dots,f_{k} \right\}$. If $f \in C(X,\mathbb{R})$ and $\mu \in M(X)$, then
$$ \overline{P}_{\mu,{\mathcal{G}} }(f) \leq P^{KP}_{\mu,\mathcal{G}}(f).$$
\end{proposition}

For the proof of the Proposition \ref{p3}, we need the following lemma (cf. \cite[Theorem 2.1]{Mat}).
\begin{lemma}\label{l2}
   (5r-lemma) Let $(X,d)$ be a compact metric space and $\mathcal{B}=\left \{ B(x_{i},r_{i}) \right \} _{i \in I}$ be a family of open (or closed) balls in $X$. Then there exists a finite or countable subfamily $\mathcal{B'}=\left \{ B(x_{i},r_{i}) \right \} _{i \in I'}$ of pairwise disjoint balls in $\mathcal{B}$ such that
   $$\bigcup_{B \in \mathcal{B} }B \subset \bigcup_{i \in I'}B(x_{i},5r_{i}) .$$
\end{lemma}

\begin{proof}[Proof of Proposition~{\upshape\ref{p3}}]
For any $s<\overline{P}_{\mu,\mathcal{G} }(Z,f)$, we can find $\varepsilon,\ \beta>0$ and $A \subseteq Z$ with $\mu(A)>0$ such that
$$\limsup_{n \to +\infty} \frac{-\log{\mu}(B_{n}(x,e^{-n\varepsilon} )) +f_{n}(x)}{| G_{n}| }>\beta +s , \ \forall x \in A.$$
Fix $\delta \in (0,\mu(A))$, we shall show that
$$P^{KP}_{\mu,\mathcal{G}}(s,2\varepsilon,\delta, f)=+\infty.$$
Let $\left\{Z_{i}\right\}_{i \in I}$ be a countable family with $\mu(\bigcup_{i}Z_{i})>1-\delta$. It follows that
$$\mu(A \cap \bigcup_{i}Z_{i}) \geq \mu(A)-\delta >0.$$
Hence there exists $i$ such that $\mu(A \cap Z_{i})>0$. For such $i$, we define
$$E_{n}=\left \{ x \in A \cap Z_{i}:\mu(B_{n}(x,e^{-n\varepsilon} ))<e^{-| G_{n} |(\beta+s)+f_{n}(x) } \right \}, \ n \in \mathbb{N}.$$
It is clear that $A \cap Z_{i}=\bigcup^{\infty}_{n=N}E_{n}$ for each $N \in \mathbb{N}$. Fix $N \in \mathbb{N}$, then $\sum_{n=N}^{\infty }\mu(E_{n}) \geq \mu(A \cap Z_{i}) $. It follows that there exists $n \geq N$ such that
$$\mu(E_{n})\geq \frac{1}{n(n+1)}\mu(A \cap Z_{i}).$$
Let $n$ be a sufficiently large integer so that $e^{n{\varepsilon}}>5$. Fix such $n$ and let 
$$\mathcal{B}=\left \{ B_{n}(x,e^{-2n\varepsilon } ):x \in E_{n} \right \}. $$  
By Lemma \ref{l2}, there exists a finite pairwise disjoint family $\left\{ B_{n}(x_{i},e^{-2n\varepsilon }) \right\}_{i \in I}$ with $x_{i} \in E_{n}$ such that
$$E_{n} \subset \bigcup_{x \in E_{n}}B_{n}(x,e^{-2n\varepsilon}) \subset \bigcup_{i \in I} B_{n}(x_{i},5e^{-2n\varepsilon}) \subset \bigcup_{i \in I} B_{n}(x_{i},e^{-n\varepsilon}).$$
Hence
$$\begin{aligned}
M^{P}_{\mathcal{G}}(N,s,2\varepsilon,Z_{i},f) &\geq M^{P}_{\mathcal{G}}(N,s,2\varepsilon,A\cap Z_{i},f) \geq M^{P}_{\mathcal{G}}(N,s,2\varepsilon,E_{n},f) \\
&\geq \sum_{i \in I}e^{-| G_{n}| s+f_{n}(x_{i})} =e^{| G_{n}|\beta}\sum_{i \in I}e^{-| G_{n}| (\beta+s)+f_{n}(x_{i})} \\
&\geq e^{| G_{n}|\beta}\sum_{i \in I}\mu(B_{n}(x_{i},e^{-n\varepsilon})) \geq e^{| G_{n}|\beta}\mu(E_{n}) \geq e^{| G_{n}|\beta} \frac{\mu(A \cap Z_{i})}{n(n+1)}.
\end{aligned}$$
By the definition of $G_{n}$, we have
$$\lim_{n \to +\infty}\frac{e^{|G_{n}|\beta}}{n(n+1)}=+\infty.$$
Letting $N \to +\infty$, we obtain that
$$M^{P}_{\mathcal{G}}(s,2\varepsilon,E,f)=\lim_{N \to +\infty} M^{P}_{\mathcal{G}}(N,s,2\varepsilon,E,f)=+\infty.$$
Then we have
$$P^{KP}_{\mu,\mathcal{G}}(s,2\varepsilon,\delta, f)=+\infty.$$
It follows that $\overline{P}_{\mu,{\mathcal{G}} }(f) \leq P^{KP}_{\mu,\mathcal{G}}(f).$
\end{proof}

Now we can prove the Theorem \ref{T4}.
\begin{proof}[Proof of Theorem~{\upshape\ref{T4}}]
Let $$ P_{\mu,\mathcal{G}}^{P}(f):=\lim_{\varepsilon  \to 0} \lim_{\delta  \to 0} \inf\left \{P^{P}_{\mathcal{G}}(\varepsilon ,Z,f):\mu(Z) \geq 1-\delta \right \}.$$
Then by the definition, we have 
$$P^{P}_{\mathcal{G}}(Z,f)\geq \sup \{P_{\mu,\mathcal{G}}^{P}(f):\mu \in M(X), \ \mu(Z)=1 \}.$$
We now demonstrate $P_{\mu,\mathcal{G}}^{P}(f) \geq P^{KP}_{\mu,\mathcal{G}}(f)$, the proof of Theorem \ref{T4} is completed by an application of Proposition \ref{p3}.
Let $P^{KP}_{\mu,\mathcal{G}}(f)>s$, then there exist $\varepsilon'>0$ and $\delta'>0$ such that
$$P^{KP}_{\mu,\mathcal{G}}(\varepsilon,\delta,f)>s, \ \forall \varepsilon \in (0,\varepsilon')\ \text{and} \ \delta \in (0,\delta').$$
This implies that $M_{\mu,\mathcal{G}}^{\mathcal{P}}(s,\varepsilon ,\delta,f)=+\infty$. For any $Z \subseteq X$ with $\mu(Z) \geq 1-\delta$, if $Z \subseteq \bigcup_{i}Z_{i}$, then $\mu(\bigcup_{i}Z_{i}) \geq 1-\delta$. Thus
$$\sum_{i=1}^{\infty}M^{P}_{\mathcal{G}}(s,\varepsilon ,Z_{i},f) \geq M^{P}_{\mathcal{G}}(s,\varepsilon ,Z,f)=+\infty.$$
Hence $M^{\mathcal{P}}_{\mathcal{G}}(s,\varepsilon ,Z,f)=+\infty$. It follows that
$$P_{\mu,\mathcal{G}}^{P}( f)=\lim_{\varepsilon  \to 0} \lim_{\delta  \to 0} \text{inf}\left \{ P^{P}_{\mathcal{G}}(\varepsilon ,Z, f):\mu(Z) \geq 1-\delta \right \} \geq s.$$
Then we have $P_{\mu,\mathcal{G}}^{P}(f) \geq P^{KP}_{\mu,\mathcal{G}}(f)$.
\end{proof}

\subsection{Upper bound}

\begin{theorem}\label{T5}
Let $G$ be a free semigroup acting on a compact metric space $(X,d)$ with finite generator set $\mathcal{G}=\left\{f_{1},f_{2},\dots,f_{k} \right\}$. If $f \in C(X,\mathbb{R})$, $Z \subset X$ is a nonempty compact set and $ P^{P}_{\mathcal{G}}(Z,f) > \|f\|$, where $\|f\|:=\sup_{x \in X}|f(x)|$, then
   $$P^{P}_{\mathcal{G}}(Z,f)\leq \sup \{ \overline{P}_{\mu,\mathcal{G} }(Z,f):\mu \in M(X), \ \mu(Z)=1 \}.$$
\end{theorem}
The following lemma is important in the proof of Theorem \ref{T5}.
\begin{lemma}\label{l3}
   Let $Z \subset X$, $\varepsilon>0$ and $s>\|f\|$. If $M^{P}_{\mathcal{G}}(s,\varepsilon,Z,f)=\infty$, then for a given finite interval $(a,b)\in [0,+\infty)$ and $n \in \mathbb{N}$, there exists a finite disjoint collection $\{\overline{B}_{n_{i}}(x_{i},e^{-n_{i}\varepsilon)}\}$ such that $x_{i}\in Z$, $n_{i} \geq n$ and $\sum_{i}e^{-s\left|G_{n_{i}}\right|+f_{n_{i}}(x_{i})}\in (a,b).$
\end{lemma}
\begin{proof}
   Let $N>n$ be large enough such that $e^{|G_{N}| (\|f\|-s)}<b-a$. Since $M^{P}_{\mathcal{G}}(s,\varepsilon,Z,f)=+\infty$, we have  $M^{P}_{\mathcal{G}}(N,s,\varepsilon,Z,f)=+\infty$. Then there exists a finite disjoint collection $\{\overline{B}_{n_{i}}(x_{i},e^{-n_{i}\varepsilon)}\}$ such that $x_{i} \in Z$, $n_{i} \geq N$ and $\sum_{i}e^{-s|G_{n_{i}}|+f_{n_{i}}(x_{i})}>b$. Since
   $$e^{-s|G_{n_{i}}|+f_{n_{i}}(x_{i})}\leq e^{-s|G_{n_{i}}|+|G_{n_{i}}| \cdot \|f\|} \leq e^{|G_{n_{i}}| (\|f\|-s)}<b-a,$$
   we can discard elements in this collection one by one until
   $$\sum_{i}e^{-s|G_{n_{i}}|+f_{n_{i}}(x_{i})} \in (a,b).$$
\end{proof}

\begin{proof}[Proof of Theorem~{\upshape\ref{T5}}]
For any $s \in (\|f\|,P^{P}_{\mathcal{G}}(Z,f))$, we have $s<P^{P}_{\mathcal{G}}(\varepsilon,Z,f)$ for all $\varepsilon \in(0,1)$. Fix $t \in (s,P^{P}_{\mathcal{G}}(\varepsilon,Z,f))$. We will construct inductively the following four sequences inspired by the work of Feng and Huang \cite{Fe}:
\begin{enumerate}[1)]
\item[1)] a sequence of finite sets $\left\{K_{i} \right\}$ with $K_{i}\subset Z$;
\item[2)] a sequence of finite measures $\left\{\mu_{i}\right\}$ with each $\mu_{i}$ being supported on $K_{i}$;
\item[3)] a sequence of positive numbers $\left\{\gamma_{i}\right\}$;
\item[4)] a sequence of integer-valued functions $\left\{m_{i}\right\}$ where $m_{i}:K_{i}\to \mathbb{N}$.
\end{enumerate}
The construction is divided into three steps:

$\textbf{Step 1.}$ Construct $K_{1}$, $\mu_{1}$, $m_{1}(\cdot )$ and $\gamma_{1}$.

Note that $M^{\mathcal{P}}_{\mathcal{G}}(t,\varepsilon,Z,f)=+\infty$. Let
$$H=\bigcup\left\{J \subset X:J \ \text{is open}, \ M^{\mathcal{P}}_{\mathcal{G}}(t,\varepsilon,Z \cap J,f)=0\right\}.$$

Then by the separability of $X$, $H$ is a countable union of the open set $J$. It implies that $M^{\mathcal{P}}_{\mathcal{G}}(t,\varepsilon,Z \cap H,f)=0$. Let $Z'=Z\setminus H=Z \cap (X\setminus H)$. We first show that for any open set $J \subset X$, either $Z' \cap J=\emptyset $ or $M^{\mathcal{P}}_{\mathcal{G}}(t,\varepsilon,Z' \cap J,f)=0$.

Suppose that $M^{\mathcal{P}}_{\mathcal{G}}(t,\varepsilon,Z' \cap J,f)=0$. Since $Z=Z' \cup (Z \cap H)$, we have
$$M^{\mathcal{P}}_{\mathcal{G}}(t,\varepsilon,Z \cap J,f) \leq M^{\mathcal{P}}_{\mathcal{G}}(t,\varepsilon,Z' \cap J,f)+M^{\mathcal{P}}_{\mathcal{G}}(t,\varepsilon,Z \cap H,f)=0.$$
Thus $J \subset H$. It follows that $Z' \cap J=\emptyset$. Since
$$M^{\mathcal{P}}_{\mathcal{G}}(t,\varepsilon,Z,f) \leq M^{\mathcal{P}}_{\mathcal{G}}(t,\varepsilon,Z \cap H,f) + M^{\mathcal{P}}_{\mathcal{G}}(t,\varepsilon,Z',f)$$
and $M^{\mathcal{P}}_{\mathcal{G}}(t,\varepsilon,Z \cap H,f)=0$, it holds that
$$M^{\mathcal{P}}_{\mathcal{G}}(t,\varepsilon,Z',f)= M^{\mathcal{P}}_{\mathcal{G}}(t,\varepsilon,Z,f)=+\infty.$$
Thus $M^{\mathcal{P}}_{\mathcal{G}}(s,\varepsilon,Z',f)=+\infty$.

Using Lemma \ref{l3}, we can find a finite set $K_{1} \subset Z'$ and an integer-valued function $m_{1}(x)$ on $K_{1}$ such that the collection $\{ \overline{B}_{m_{1}(x)}(x,e^{-m_{1}(x)\varepsilon}) \}_{x \in K_{1}}$ is disjoint and
$$\sum_{x \in K_{1}}e^{|G_{m_{1}(x)}|s+f_{m_{1}(x)}(x)} \in (1,2).$$
Define
$$\mu_{1}=\sum_{x \in K_{1}}e^{|G_{m_{1}(x)}|s+f_{m_{1}(x)}(x)} \delta_{x},$$
where $\delta_{x}$ denotes the Dirac measure at $x$. Take $\gamma_{1}>0$ small enough such that for any $z \in \overline{B}(x,\gamma_{1})$, we have
\begin{equation}
(\overline{B}(z,\gamma_{1}) \cup  \overline{B}_{m_{1}(x)}(z,e^{-m_{1}(x)\varepsilon})) \cap (\bigcup_{y \in K_{1} \setminus \{x\}}\overline{B}(y,\gamma_{1}) \cup  \overline{B}_{m_{1}(y)}(y,e^{-m_{1}(y)\varepsilon}))=\emptyset.\label{eq:4.3}
\end{equation}
Since $K_{1} \subset Z'$, we have
$$M^{\mathcal{P}}_{\mathcal{G}}(t,\varepsilon,Z \cap B(x,\gamma_{1}/4),f) \geq M^{\mathcal{P}}_{\mathcal{G}}(t,\varepsilon,Z' \cap B(x,\gamma_{1}/4),f)>0,$$
for any $x \in K_{1}$.

$\textbf{Step 2.}$ Construct $K_{2}$, $\mu_{2}$, $m_{2}(\cdot )$ and $\gamma_{2}$.

By \eqref{eq:4.3}, the elements of the family of balls $\{ \overline{B}(x,\gamma_{1})\}_{x \in K_{1}}$ are pairwise disjoint. For each $x \in K_{1}$, since $M^{\mathcal{P}}_{\mathcal{G}}(t,\varepsilon,Z \cap B(x,\gamma_{1}/4),f)>0$, we can construct a finite set $E_{2}(x) \subset Z \cap B(x,\gamma_{1}/4)$ and integer-valued function
$$m_{2}(x):E_{2}(x) \to \mathbb{N} \cap [\text{max}\{m_{1}(y):y \in K_{1}\},+\infty)$$
such that
\begin{enumerate}
\item[(a)] $M^{\mathcal{P}}_{\mathcal{G}}(t,\varepsilon,Z \cap J,f)>0$ for each open set $J$ with $J \cap E_{2}(x)\neq \emptyset$;
\item[(b)] The elements in $\left \{ B_{m_{2}(y)}(y,e^{-m_{2}(y)\varepsilon} ) \right \}_{y \in E_{2}(x)} $ are disjoint and
$$ \mu_{1}(\{x\})<\sum_{y \in E_{2}(x)}e^{-|G_{m_{2}(y)}|s+f_{m_{2}(y)}(y)}<(1+2^{-2})\mu_{1}(\{x\}). $$
\end{enumerate}
To see it, we fix $x \in K_{1}$ and denote $F=Z \cap B(x,\gamma_{1}/4)$. Let
$$H_{x}=\bigcup \left\{ J \subset X :J \ \text{is open},\ M^{\mathcal{P}}_{\mathcal{G}}(t,\varepsilon,F \cap J,f)=0 \right\}.$$
Set $F'=F \setminus H_{x}$. Then by Step 1, we can show that
$$M^{\mathcal{P}}_{\mathcal{G}}(t,\varepsilon,F',f)=M^{\mathcal{P}}_{\mathcal{G}}(t,\varepsilon,F,f)>0$$
and $M^{\mathcal{P}}_{\mathcal{G}}(t,\varepsilon,F' \cap J,f)>0$ for any open set $J$ with $J \cap F' \neq \emptyset$. Since $s<t$, we have $M^{\mathcal{P}}_{\mathcal{G}}(s,\varepsilon,F',f)=+\infty$.

Using Lemma \ref{l3} again, we can find a finite set $E_{2} \subset F'$ and a map
$$m_{2}(x):E_{2}(x) \to \mathbb{N} \cap [\text{max}\{m_{1}(y):y \in K_{1}\},+\infty)$$
so that (b) holds. If $J \cap E_{2}(x) \neq \emptyset$ and $J$ is an open set, then $J \cap F' \neq \emptyset$. Hence
$$M^{\mathcal{P}}_{\mathcal{G}}(t,\varepsilon,Z \cap J,f) \geq M^{\mathcal{P}}_{\mathcal{G}}(t,\varepsilon,F'\cap J,f)>0.$$
Thus (a) holds.
Since the elements of the family $\{ \overline{B}(x,\gamma_{1})\}_{x \in K_{1}}$ are pairwise disjoint, $E_{2}(x) \cap E_{2}(y)=\emptyset$ for different points $x,y \in K_{1}$. Define $K_{2}=\bigcup_{x \in K_{1}}E_{2}(x)$ and
$$\mu_{2}=\sum_{x \in K_{2}}e^{|G_{m_{2}(x)}|s+f_{m_{2}(x)}(x)} \delta_{x}.$$
By \eqref{eq:4.3} and (b), the elements in $\{ \overline{B}_{m_{2}(x)}(x,e^{-m_{2}(x)\varepsilon}) \}_{x \in K_{2}}$ are pairwise disjoint. Hence we can take $\gamma_{2} \in (0,\gamma_{1}/4)$ small enough such that for any  $x \in K_{2}$ and the function $z:K_{2} \to X$ with $d(x,z(x))\leq \gamma_{2}$, we have
$$(\overline{B}(z(x),\gamma_{2}) \cup  \overline{B}_{m_{2}(x)}(z(x),e^{-m_{2}(x)\varepsilon})) \cap (\bigcup_{y \in K_{2} \setminus \{x\}}\overline{B}(z(y),\gamma_{2}) \cup  \overline{B}_{m_{2}(y)}(z(y),e^{-m_{2}(y)\varepsilon}))=\emptyset.$$
Since $x \in K_{2}$, there exists $y \in K_{1}$ such that $x \in E_{2}(y)$. By (a), we have
$$M^{\mathcal{P}}_{\mathcal{G}}(t,\varepsilon,Z \cap B(x,\gamma_{2}/4),f)>0,$$
for each $x \in K_{2}$.

$\textbf{Step 3.}$ Assume that $K_{i}$, $\mu_{i}$, $m_{i}(\cdot )$ and $\gamma_{i}$ have been constructed, for $i=1,2,\cdots,p$. In particular, suppose that for any $x \in K_{p}$ and the function $z:K_{p} \to X$ with $d(x,z(x))<\gamma_{p}$,
\begin{equation}
(\overline{B}(z(x),\gamma_{p}) \cup  \overline{B}_{m_{p}(x)}(z(x),e^{-m_{p}(x)\varepsilon})) \cap (\bigcup_{y \in K_{p} \setminus \{x\}}\overline{B}(z(y),\gamma_{p}) \cup  \overline{B}_{m_{p}(y)}(z(y),e^{-m_{p}(y)\varepsilon}))=\emptyset \label{eq:4.4}
\end{equation}
and $M^{\mathcal{P}}_{\mathcal{G}}(t,\varepsilon,Z \cap B(x,\gamma_{p}/4),f)>0$. We shall construct $K_{p+1}$, $\mu_{p+1}$, $m_{p+1}(\cdot )$ and $\gamma_{p+1}$ in a way being similar to Step 2.

Note that the elements in $\{ \overline{B}(x,\gamma_{p}) \}_{x \in K_{p}}$ are pairwise disjoint. Since
$$M^{\mathcal{P}}_{\mathcal{G}}(t,\varepsilon,Z \cap B(x,\gamma_{p}/4),f)>0,$$
for each $x \in K_{p}$,  we can construct as in Step 2, a finite set
$$E_{p+1}(x) \subset Z \cap B(x,\gamma_{p}/4)$$
and an integer-valued function
$$m_{p+1}(x):E_{p+1}(x) \to \mathbb{N} \cap [\text{max}\{m_{p}(y):y \in K_{p}\},+\infty),$$
such that
\begin{itemize}
\item[(c)] $M^{\mathcal{P}}_{\mathcal{G}}(t,\varepsilon,Z \cap J,f)>0$, for each open set $J$ with $J \cap E_{p+1}(x)\neq \emptyset$.
\item[(d)] The elements in $\left \{ B_{m_{p+1}(y)}(y,e^{-m_{p+1}(y)\varepsilon} ) \right \}_{y \in E_{p+1}(x)} $ are disjoint, and
$$ \mu_{p}(\{x\})<\sum_{y \in E_{p+1}(x)}e^{-|G_{m_{p+1}(y)}|s+f_{m_{p+1}(y)}(y)}<(1+2^{-p-1})\mu_{p}(\{x\}). $$
\end{itemize}
It is easy to see that $E_{p+1}(x) \cap E_{p+1}(y)=\emptyset$ for any $x,y \in K_{p}$ with $x \neq y$. Define
$$K_{p+1}=\bigcup_{x \in K_{p}}E_{p+1}(x)$$ and
$$\mu_{p+1}=\sum_{x \in K_{p+1}}e^{|G_{m_{p+1}(x)}|s+f_{m_{p+1}(x)}(x)} \delta_{x}.$$
By \eqref{eq:4.4} and (d), the elements in $\{ \overline{B}_{m_{p+1}(x)}(x,e^{-m_{p+1}(x)\varepsilon}) \}_{x \in K_{p+1}}$ are pairwise disjoint. Hence we can take $\gamma_{p+1} \in (0,\gamma_{p}/4)$ small enough such that for any  $x \in K_{p+1}$ and the function $z:K_{p+1} \to X$ with $d(x,z(x))\leq \gamma_{p+1}$, we have
\begin{align}
\nonumber&(\overline{B}(z(x),\gamma_{p+1}) \cup  \overline{B}_{m_{p+1}(x)}(z(x),e^{-m_{p+1}(x)\varepsilon})) \cap \\\nonumber&(\bigcup_{y \in K_{p+1} \setminus \{x\}}\overline{B}(z(y),\gamma_{p+1}) \cup  \overline{B}_{m_{p+1}(y)}(z(y),e^{-m_{p+1}(y)\varepsilon}))=\emptyset.
\end{align}
Since $x \in K_{p+1}$, there exists $y \in K_{p}$ such that $x \in E_{p+1}(y)$. Thus by (c),
$$M^{\mathcal{P}}_{\mathcal{G}}(t,\varepsilon,Z \cap B(x,\gamma_{p+1}/4),f)>0$$
for each $x \in K_{p+1}$.

We summarize their properties as follow:
\begin{itemize}
\item[(e)] For each $i$, the family $\mathcal{F}_{i}=\{\overline{B}(x,\gamma_{i}):x \in K_{i}\}$ is disjoint. For every $B \in \mathcal{F}_{i+1}$, there exists $x \in K_{i}$ such that $B \subset \overline{B}(x,\gamma_{i}/2);$
\item[(f)] For each $x \in K_{i}$ and $x' \in \overline{B}(x,\gamma_{i})$,
\begin{equation}
\overline{B}_{m_{i}(x)}(x',e^{-m_{i}(x)\varepsilon})\cap \bigcup_{y \in K_{i} \setminus \{x\}}\overline{B}(y,\gamma_{i})= \emptyset \label{eq:4.5}
\end{equation}
and
\begin{equation}
\begin{aligned}
\mu_{i}(\overline{B}(x,\gamma_{i}))=e^{-|G_{m_{i}(x)}|s+f_{m_{i}(x)}(x)}& \leq \sum_{y \in E_{i+1}(x)}e^{-|G_{m_{i+1}(y)}|s+f_{m_{i+1}(y)}(y)}\\
&\leq (1+2^{-i-1})\mu_{i}(\overline{B}(x,\gamma_{i})), \label{eq:4.6}
\end{aligned}
\end{equation}
where $E_{i+1}(x)=\overline{B}(x,\gamma_{i})\cap K_{i+1}$.
\end{itemize}
   By \eqref{eq:4.6}, for any $F_{i} \in \mathcal{F}_{i}$, we have
   $$\begin{aligned}
      \mu_{i}(F_{i}) \leq \mu_{i+1}(F_{i})&=\sum_{F \in \mathcal{F}_{i+1}:F \subset F_{i}}\mu_{i+1}(F)\\
      &\leq \sum_{F \in \mathcal{F}_{i+1}:F \subset F_{i}}(1+2^{-i-1})\mu_{i}(F)\\
      &=(1+2^{-i-1})\sum_{F \in \mathcal{F}_{i+1}:F \subset F_{i}}\mu_{i}(F)\\
      &\leq (1+2^{-i-1})\mu_{i}(F_{i}).
   \end{aligned}$$
   Using the above inequality repeatedly, we have for any $j>i$,
   \begin{equation}
      \mu_{i}(F_{i}) \leq \mu_{j}(F_{i}) \leq \prod_{n=i+1}^{j}(1+2^{-n})\mu_{i}(F_{i})\leq C\mu_{i}(F_{i}), \ \forall F_{i}\in \mathcal{F}_{i}, \label{eq:4.7}
   \end{equation}
   where $C=\prod_{n=1}^{+\infty}(1+2^{-n})<+\infty$.

   Let $\hat{\mu}$ be a limit point of $\{\mu_{i}\}$ in a weak$^{*}$ topology. Let
   $$K^{*}=\bigcap_{n=1}^{\infty}\overline{\bigcup_{i\geq n}K_{i}}=\lim_{n \to +\infty}\overline{\bigcup_{i\geq n}K_{i}}.$$
   Then $\hat{\mu}$ is supported on $K^{*}$, $K^{*}\subset Z$ and for any $i \in \mathbb{N}$, $K^{*}\subset \bigcup_{x \in K_{i}}B(x,\gamma_{i})$. By \eqref{eq:4.7},
   $$\begin{aligned}
      e^{-| G_{m_{i}(x)}|s+f_{m_{i}(x)}(x) }&=\mu_{i}(\overline{B}(x,\gamma_{i})) \leq \hat{\mu}(\overline{B}(x,\gamma_{i}))\\
      &\leq C\mu_{i}(\overline{B}(x,\gamma_{i})) =C e^{-| G_{m_{i}(x)}|s+f_{m_{i}(x)}(x) }, \ \forall x \in K_{i}.
   \end{aligned}$$
   In particular,
   $$1 \leq \sum_{x \in K_{1}}\mu_{1}(B(x,\gamma_{1}))\leq \sum_{x \in K_{1}}\hat{\mu}(B(x,\gamma_{1}))=\hat{\mu}(K^{*})\leq \sum_{x \in K_{1}}C\mu_{1}(B(x,\gamma_{1}))\leq 2C.$$
   By \eqref{eq:4.5}, for every $x \in K_{i}$ and $x'\in \overline{B}(x,\gamma_{i})$,
   $$\hat{\mu} (\overline{B}_{m_{i}(x)}(x',e^{-m_{i}(x)\varepsilon})) \leq \hat{\mu}(\overline{B}(x,\gamma_{i})) \leq Ce^{-|G_{m_{i}(x)}|s+f_{m_{i}(x)}(x) }.$$
   For each $x' \in K^{*}$ and $i \in \mathbb{N}$, there exists $x\in K_{i}$ such that $x'\in \overline{B}(x,\gamma_{i})$. Thus
   $$\hat{\mu} (\overline{B}_{m_{i}(x)}(x',e^{-m_{i}(x)\varepsilon})) \leq C e^{-| G_{m_{i}(x)}|s+f_{m_{i}(x)}(x) }.$$
   Let $\mu=\hat{\mu}/\hat{\mu}(K^{*})$. Then $\mu \in M(x)$ and $\mu(K^{*})=1$. Moreover, for all $x' \in K^{*}$, there exists a sequence $\{k_{i}\}_{i\geq 1}$ with $k_{i} \to +\infty$ such that
   $$\mu(B_{k_{i}}(x',e^{-k_{i}\varepsilon})) \leq \frac{Ce^{-| G_{k_{i}}|s+f_{k_{i}}(x') }}{\hat{\mu}(K^{*})}.$$
   It implies that
   $$\frac{-\log{\frac{\hat{\mu}(K^{*}) }{C} } -\log{\mu(B_{k_{i}}(x',e^{-k_{i}\varepsilon}))}+f_{k_{i}}(x')}{| G_{k_{i}}| } \geq s.$$
   Letting $k_{i} \to +\infty$, we get $\overline{P}_{\mu,\mathcal{G}}(\varepsilon,Z,f) \geq s$.
 \end{proof}

\vskip 0.3cm {\bf Acknowledgements}\quad The authors would like to thank the anonymous referees for their valuable comments and suggestions. The first author is supported by  NNSF of China (Grant No.12201120).





\end{sloppypar}
\end{document}